\documentclass[letterpaper, 10 pt, conference]{ieeeconf}  

\IEEEoverridecommandlockouts                              

\usepackage{graphics} 
\usepackage{epsfig} 
\usepackage{amsmath} 
\usepackage{amsmath} 
\usepackage{amssymb}  
\usepackage{cite}
\usepackage{fullpage}
\usepackage{fancyhdr}
\usepackage{color, subfigure,multirow}

\newtheorem{defn}{Definition}

\newtheorem{lem}{Lemma}

\newtheorem{prop}{Proposition}

\newenvironment{sketch-proof}{{\it Sketch of the Proof:\ }}{ \hfill}

\newcommand{\R}{{\mathbb R}}

\renewcommand{\Pr}{{\mathbb{P}}}

\newcommand{\diag}{{\text{diag}}}
\newcommand{\blkdiag}{{\text{diag}}}

\newcommand{\figref}[1]{Figure~\ref{#1}}
\newcommand{\figrefs}[2]{Figures~\ref{#1},\ref{#2}}

\title{On Projection-Based Model Reduction of Biochemical Networks\\
 Part I: The Deterministic Case}
\author{Aivar Sootla$^{1}$ and James Anderson$^{2}$
\thanks{$^{1}$AS is with the Department of Bioengineering, Imperial College London, UK {\tt\small a.sootla@imperial.ac.uk}}%
\thanks{$^{2}$JA is with St. John's College and the Department of Engineering Science, University of Oxford, Parks Road, OX1 3JP {\tt\small james.anderson@eng.ox.ac.uk}}%
}

\begin{document}
\maketitle
\thispagestyle{empty}
\pagestyle{empty}
\begin{abstract}
This paper addresses the problem of model reduction for dynamical system models that describe biochemical reaction networks. Inherent in such models are properties such as stability, positivity and network structure. Ideally these properties should be preserved by model reduction procedures, although traditional projection based approaches struggle to do this. We propose a projection based model reduction algorithm which uses generalised block diagonal Gramians to preserve structure and positivity. Two algorithms are presented, one provides more accurate reduced order models, the second provides easier to simulate reduced order models. The results are illustrated through numerical examples.
\end{abstract}

\section{Introduction}
Biochemical reaction networks are most appropriately modelled as stochastic systems. Typically they take the form of an infinite dimensional, continuous time Markov Chain which describes the time evolution of a probability density function of the concentration of the reactants. The Chemical Master Equation (CME)
\begin{equation*}
\frac{\partial \Pr(n,t)}{\partial t }=\Omega \sum_{i=1}^R \hat f(n-S_i,\Omega)-\hat f(n,\Omega))\Pr(n,t)
\end{equation*}
describes how a reaction network composed of: $R$ reactions, in a compartment of volume $\Omega$ with a stoichiometry matrix $S$ ($S_i$ denoting the $i^{\text{th}}$ column); $\hat f$ the flux vector; $n$ the vector containing the number of molecules $n_i$ of species $i$; and $\Pr(n,t)$ is the probability of the vector of molecules $n$ at time $t$ changes with time.

For reaction networks with just a few species even simulating the CME can be intractable. This paper is a first step towards an automated procedure to compute efficient reduced order models for stochastic biochemical models. It is assumed that the starting point for the algorithms presented here is a nonlinear, possibly high dimensional, but deterministic dynamical system. Part 2 of this paper\footnote{This paper is completely self-contained and does not require any of the material from Part II.} \cite{SootlaAndACC} and the references therein describes how and under what assumptions one can approximate the CME by a deterministic dynamical system. 

The focus of this paper is to describe a projection based algorithm for reducing the state dimension of a dynamical system while preserving certain desirable features such as stability, positivity and network structure. Standard model reduction techniques make use of the fact that frequently states (species concentrations) evolve over multiple time scales \cite{Tikhonov1952,kokotovic1987singular}. The basic idea with such methods is to treat the \emph{fast} states as being at steady state, thus obtaining an algebraic expression which can be then substituted into the slow state dynamics. Such approaches are referred to as time scale separation or quasi-steady state assumption (QSSA) methods.

More common in the control literature is the use of \emph{projection} based model order reduction \cite{Moore1981,Scherpen93,nonlinearkrylov,AntoulasBook}. Balancing-based projection methods follow a two step procedure; first a state-space transformation is found which aligns the controllability and observability ellipsoids, then the states which are least controllable and observable are truncated yielding a reduced order model. In some cases, provided the initial full order model was stable it can be shown that the reduced model is stable too. It is often possible to \emph{a priori} determine the error bound in an appropriate choice of norm between the full and reduced model. The major drawback of projection approaches is that the states in the transformed coordinate system are linear combinations of all the other states, thus the physical meaning of a state is lost. Recently structure preserving reduction algorithms have been proposed based upon coprime factorisation~\cite{LiPaganini}, structured Gramians~\cite{Sandberg09}, $\mathcal H_{\infty}$ optimisation~\cite{Sootla2012structured,hinfstruct} and novel energy functions~\cite{Sootla2012positive} that attempt to avoid such problems.  The work in this paper most closely resembles the spirit of~\cite{Sandberg09}, however in addition to preserving network structure would like to preserve the monotonicity, when possible. 
\section{Problem Statement}
The standard model reduction problem takes the following form: given a stable dynamical system
\begin{equation}\label{eq:nl}
  \begin{aligned}
    \dot x&= f(x, u) \\
         y&= h(x)
  \end{aligned}
\end{equation}
where $x\in \R^n$ is the state vector, and the equilibrium point of interest is without loss of generality $x_{ss}=0_{n\times 1}$.  Construct a dynamical system 
\begin{equation}\label{eq:nl_red}
  \begin{aligned}
    \dot{\tilde{x}}&= \tilde f(\tilde{x},u) \\
       \tilde y&= \tilde h(\tilde x) 
  \end{aligned} 
\end{equation}
where $\tilde{x}\in \R^k$ with $k<n$ and the \emph{error} between \eqref{eq:nl}--\eqref{eq:nl_red} is small in some appropriate norm. When $f$ is nonlinear the reduction problem is in general intractable, see \cite{Scherpen93} for nonlinear input-affine balancing and \cite{nonlinearkrylov} for SISO nonlinear moment matching approaches. In this paper we shall deal with linearisations about a given operating point and input and adapt classical methods (cf. \cite{Glover1984}) to preserve desirable system properties as outlined in the next section. In order to simplify some derivations, we assume that $h(x) = C x$, where $C$ is a constant matrix.
											
\subsection{Structured Projectors}
We approximate the system \eqref{eq:nl} around the stable steady-state $x_{ss}$ with a constant control signal $u_{ss}$. 

Consider a system 
\begin{equation}
  \label{eq:sys}
  \begin{aligned}
    \dot x&=A x+B u\\
         y&=C x
  \end{aligned}
\end{equation}
where the drift matrix $A $ and  input map $B$ are given by
\begin{equation*}
A=\dfrac{\partial f(x,u)}{\partial x}\Bigl|_{x = x_{s s}, u = u_{s s}}, \text{ } B =  \dfrac{\partial f(x,u)}{\partial u}\Bigl|_{x = x_{s s}, u = u_{s s}} .
\end{equation*}
Note that $A $ is Hurwitz by assumption. The linearised system \eqref{eq:sys} can then be partitioned as follows:
\begin{equation} \small
  \label{eq:part}
x = \begin{pmatrix}  x_1 \\ x_2 \end{pmatrix}~ 
A =\begin{pmatrix} A_{11} & A_{12} \\ A_{21} & A_{22} \end{pmatrix}~ 
B = \begin{pmatrix} B_{1} \\ B_{2} \end{pmatrix}~ 
C^T = \begin{pmatrix}  C_{1}^T \\ C_{2}^T \end{pmatrix},
\end{equation}
where $x_1 \in {\R}^{n-k}$, $x_2 \in {\R}^{k}$, and the matrices $A$, $B$ and $C$ are partitioned conformally. The next step is to computed structured Gramians, which are obtained as solutions to Lyapunov inequalities
\begin{equation}
  \label{eq:lyap_ineq}
  \begin{aligned}
    A P + P A^T + B B^T&\le 0 \\
    Q A + A^T Q + C^T C&\le 0
  \end{aligned}
\end{equation}
with $P\ge 0$, $Q\ge 0$, subject to the same partitioning as the states:
\begin{equation}
\label{eq:gram_part}
P =\begin{pmatrix} P_{11} & 0_{n-k,k}  \\ 0_{k, n-k}  & P_{22} \end{pmatrix}, \text{ }
Q =\begin{pmatrix} Q_{11} & 0_{n-k,k}  \\ 0_{k, n-k}  & Q_{22} \end{pmatrix}.
\end{equation}
In the following section, we make a case why constraining the generalised Gramians $P$ and $Q$ to be block diagonal is not a restrictive assumption for biochemical networks.

If the states $x_2$ are to be approximated, the transformation $T$ is composed as follows:
\begin{equation}
  \label{eq:sst}
     T = \begin{pmatrix} I_{n-k} & 0_{n-k,k} \\ 0_{n-k,k} & T_{22} \end{pmatrix},
\end{equation}
where $T_{22}$ is such that 
\[
T_{22}^{-1} P_{22} T_{22}^{-T}=T_{22}^TQ_{22} T_{22} =\Sigma_{22},
\]
where $\Sigma_{22}$ is diagonal. According to standard tools~\cite{Sandberg09}, we choose the states to truncate according to the magnitude of the values of the diagonal of $\Sigma_{22}$. Assume $r$ states are to be reduced, let $W_{2 2}$ be the first $k-r$ columns of $T_{22}$, while $W_{2 2}^r$ are the rest $r$ columns of $T_{2 2}$, let also $V_{2 2}$ be the first $k-r$ columns of $T_{22}^{-1}$, while $W_{2 2}^r$ are the rest $r$ columns of $T_{2 2}^{-1}$. Now, the projectors can be obtained as follows
{ \begin{equation}
  \label{eq:proj}\small
  \begin{gathered}
   W = \begin{pmatrix} I_{n-k} & 0_{n-k,k-r} \\ 0_{k-r, n-k} & W_{22} \end{pmatrix} \quad   W_r = \begin{pmatrix}  0_{n-k,r} \\ W_{22}^r \end{pmatrix} \\
   V = \begin{pmatrix} I_{n-k} & 0_{n-k,k-r} \\ 0_{k-r, n-k}  & V_{22} \end{pmatrix} \quad   V_r = \begin{pmatrix}  0_{n-k,r} \\ V_{22}^r \end{pmatrix}
  \end{gathered}
\end{equation}}

\subsection{Computing a transformation for networks with monotone dynamics}
It is assumed that the dynamics of the biochemical network models we are interested in can be captured via a stoichiometric matrix $S\in\R^{n \times m}$ and flux vector $f(x)\in \R^{m\times 1}$, where $n$ is the number of species, $m$ the number of reactions that take place and $x$ the vector of species concentrations. The uncontrolled system then takes the form $\dot x = S f(x)$. We limit our focus to systems with cooperative or monotone with respect to the positive orthant $\R_{\ge 0}^n$ dynamics. This means that the stoichiometry matrix $S$ and the fluxes $f(x)$ form a \emph{cooperative} dynamical system. The following definitions make the preceding comments precise.
\begin{defn}
Consider the dynamical system $\dot{x}=r(x)$ where $r$ is locally Lipschitz, $r:\R^n_{\ge 0} \rightarrow \R^n$ and $r(0)=0$. The associated flow map is $\rho: \R_{\ge 0} \times \R^n_{\ge 0} \rightarrow \R^n$. The system is said to be monotone (w.r.t. $\R^n_{\ge 0}$) if $x\le y\Rightarrow \rho(t,x)\le \rho(t,y)$ for all $t\in \R_{\ge 0}$. 
\end{defn}

\begin{defn}
A matrix $M\in \R^{n\times n}=\left\{m_{ij}\right\}$ is said to be Metzler if $m_{ij}\ge 0$ for all $i\neq j$.
\end{defn}

The following proposition is a simplified reformulation of a known result (cf. \cite{smith2008monotone}), which establishes a straightforward test for cooperativity:
\begin{prop}
 A system $\dot x = r(x)$ is monotone with respect to the positive orthant if and only if 
\[
\frac{\partial(r^{i}(x))}{\partial x_j} \ge 0 \quad \forall i \ne j \forall x
\]
Or simply put, the Jacobian of $r(x)$ is a Metzler matrix for all $x$ in $\R_{\ge 0}^n$. 
\end{prop}

A generalisation  can be defined with respect to any orthant by mapping this orthant onto the positive one by 
a linear transformation $P: \R^n \rightarrow \R^n$, where $P = \diag((-1)^{\varepsilon_1}, \dots, (-1)^{\varepsilon_n})$ for some $\varepsilon_i$. 

Hence, after linearisation around a steady-state we have system \eqref{eq:sys} with additional constraint that the drift matrix $A$ is Metzler. We do not require $B$ and $C$ to be nonnegative matrices as is typically the case when studying positive systems.  

A simple example illustrates the concept. Consider the system
\begin{eqnarray*}
\dot{x}_1 = -x_2^3, && \dot{x}_2 = \frac{1}{1+x_1},\\
\dot{x}_3 = x_1+2x_2,&&\dot{x}_4 = x_4(\alpha x_3-\beta x_2),
\end{eqnarray*}
with $\alpha, \beta >0$. The non-zero, off diagonal elements of the Jacobian are $-3x_2^2$, $1,2$, $-\beta x_4$, $\alpha x_4$ and $\alpha x_3$. None of which change sign for $x_i>0$ and thus can be mapped to the positive othant.

In this setting our model reduction problem is formulated as replacing the states $x_2$ with a single state, while preserving stability and the Metzler property of the drift matrix. In order to obtain the reduced order model, the generalised Lyapunov equations with block-diagonal Gramians are employed. Hence the first task is to ensure the existence of such generalised Gramians.
\begin{lem}\label{lem:exist}
Consider the system \eqref{eq:sys} with an asymptotically stable, Metzler drift matrix. Let the system \eqref{eq:sys} be partioned as in \eqref{eq:part}.  Let $P$, $Q$ be generalised Gramians satisfying the Lyapunov inequalities \eqref{eq:lyap_ineq}.
Then there always exist nonnegative and nonnegative semidefinite matrices $P$ and $Q$ satisfying \eqref{eq:lyap_ineq} and the partionining as in \eqref{eq:gram_part}.
\end{lem}
\begin{proof}
See appendix.
\end{proof}
Given these properties we are ready to produce a model reduction algorithm:\\
\textbf{\underline{Reduction Algorithm:}}
\begin{enumerate}
\item Solve \eqref{eq:lyap_ineq} and obtain the matrices $P$ and $Q$ with the structure described by \eqref{eq:gram_part}.
\item Compute a balancing transformation $T_{22}$ for matrices $P_{22}$ and $Q_{22}$ as in \eqref{eq:sst}.
\item Define the projectors $W$ and $V$ as in \eqref{eq:proj} for $r$ equal to one, let $w = W_{22}$ and $v = V_{22}$.
\item Compute the matrices of the truncated reduced order model as follows
\[
\begin{aligned}
 A_t = V^T A W, \quad B_t = V^T B, \quad C_t = C W.  \\
\end{aligned}
\]
\end{enumerate}
The exact expression for matrices $A_t$, $B_t$ and $C_t$ are
\begin{equation}\label{eq:romss}
  \begin{gathered}
    A_t = \begin{pmatrix} A_{11} & A_{12} w \\ v^T A_{21} & v^T A_{22} w \end{pmatrix}\\ B_t = \begin{pmatrix} B_{1} \\ v^T B_{2} \end{pmatrix} \quad \quad
C_t^T = \begin{pmatrix}  C_{1}^T \\ w^T C_{2}^T \end{pmatrix}.
  \end{gathered}
\end{equation}
\begin{lem} \label{lem:main}
Let $P$ and $Q$ be block-diagonal matrices satisfying the conditions of Lemma~\ref{lem:exist}. Assume the matrix $P_{22} Q_{22}$ is irreducible. Let $T_{22}$ be a transformation such that $T_{22}^{-1} P_{22} T_{22}^{-T}=T_{22}^TQ_{22} T_{22}=\Sigma$, where $\Sigma$ is diagonal and $\Sigma_{1 1} \ge \Sigma_{2 2} \ge \dots \ge \Sigma_{k k}$. Let $w$ be the first column of $T$ and $v$ be the first column of $T_{22}^{-T}$. There exist such a balancing transformation $T_{2 2}$ that
\begin{enumerate}
\item The vectors $w$ and $v$ are nonnegative. 

\item The matrix $A_t$ from \eqref{eq:romss}
    is stable and Metzler.
\item Let $G$ be the full order model with a state-space realisation $A$, $B$, $C$ and $G_r$ be the reduced order model with a state-space realisation $A_t$, $B_t$, $C_t$ defined in \eqref{eq:romss}. Then 
\[
\|G-G_r\|_{\infty} \le 2 \sum\limits_{i = 2}^{k} \Sigma_{i i}
\]
\end{enumerate}
\end{lem}
\begin{proof}
See appendix.
\end{proof}

\subsection{Approximation Procedures}
The idea of using generalised structured Gramians for structured reduction is not new, in \cite{Sandberg09} they are used in an LFT framework for example. Moreover, the class of models, for which block-diagonal Gramians exists, is not rich and not many necessary conditions for the existence of block-diagonal Gramians are known. However, in the context of biochemical networks for the systems with monotone dynamics such Gramians always exist (as was just shown). Moreover, there are biochemical networks, which are not monotone, but block-diagonal Gramians still exist. Indeed it is believed that many biochemical reaction networks which do not posses the monotonicity property are actually \emph{near monotone} \cite{sontag2007monotone}. We provide an example of a non-monotone system which does admit a block diagonal Gramian in Section \ref{ex:gly}.

Define the transformed variable as $z=Tx$ and denote by $z_r$ be the species to be removed from the model, and $z_s$ the states of the reduced order model. Now the equations approximating the full order dynamics 
 can be computed as follows:
\begin{equation}
  \label{eq:red-met-1}
   \begin{aligned}
   \dot z_m& = V^Tf(W z_m + W_r z_r, u)  \\
   \dot z_r & = V_r^Tf(W z_m + W_r z_r, u)=0 \\
   y_r^d& = \Omega C (W z_m + W_r z_r)\\
 \end{aligned}
\end{equation}
Computing the root $z_r$ satisfying the algebraic-differential equation can be a computationally expensive task. Moreover, introducing the algebraic constraints may result in a stiff system, which are hard to simulate. Therefore, we also propose a truncation method, where we assume that the system is near the steady-state $z_r^0$:
\begin{equation}
  \label{eq:red-met-2}
   \begin{aligned}
   \dot z_m&= V^T f(W z_m + W_r z_r^0, u)  \\
   y_r^d&=   C (W z_m + W_r z_r^0) 
 \end{aligned} 
\end{equation}
For future reference we will refer to \eqref{eq:red-met-1} as the \emph{reduction method} and \eqref{eq:red-met-2} as the \emph{truncation method}. 

Observe that preservation of the (global) monotonicity of the reduced nonlinear system is probably not possible in general using static space-space transformations. Consider the dynamics of $z_m$ in \eqref{eq:red-met-1} with $u=0$. Let $\Gamma = W z_m + W_r z_r$. By assumption we have that $f$ is monotonic. In order for the unforced system in \eqref{eq:red-met-1} to be monotonic it needs to be shown that the Jacobian of
\begin{equation*}
\dot z_m = V_{22}^Tf_1(\Gamma,0)
\end{equation*}
given by
\begin{equation}\label{eq:Jf1}
\frac{\partial V_{22}^Tf_1^i(\Gamma(z_m),0) }{\partial (z_m)_j}\ge 0 \quad\forall i\neq j \text{ and } \forall z_m
\end{equation}
where the vector field $f$ is partitioned into $[f_1^T,f_2^T]^T$ conformally with $[z_m^T,z_r^T]$. Even for the simple case of $k-r=1$ where $V_{22}$ is simply the first column of $T_{22}^{-1}$ it is difficult to determine the underlying assumptions one would need to impose on $f$ to ensure \eqref{eq:Jf1} is satisfied.

\section{Examples}
\subsection{A Cautionary Toy Example: Topology Matters\label{ex:top}}
The first network we consider consists of four species, see \figref{fig:toy-ex}. One can interpret the species $S_1$ and $S_3$ as mRNA, and $S_2$ and $S_4$ as the corresponding proteins. As a consequence, we set the degradation rates of species $S_1$ and $S_3$ to be larger than the degradation rates of species $S_2$ and $S_4$. This however, does not necessarily imply the dynamics of $S_1$ and $S_3$ are changing on a faster time-scale than the dynamics of $S_2$ and $S_4$. We apply the standard time-scale separation technique to the network as well as the proposed reduction method \eqref{eq:red-met-1} with different ad-hoc partitions of the states\footnote{Determining \emph{a priori} appropriate partitions of a dynamical system is an open research question, see \cite{anderson2012decomposition} for example.}: lump together species $S_1$ and $S_3$ (see, \figref{fig:toy-ex-1}), lump together $S_1$ and $S_2$ (see, \figref{fig:toy-ex-2}), and finally, lump together $S_1$ and $S_2$, and simultaneously lump together $S_3$ and $S_4$ (see, \figref{fig:toy-ex-3}). The purpose of this example is signify the importance of an appropriate partitioning. The model of the network is as follows:
\begin{figure}[t]
   \centering
  \subfigure[The full order model]{\includegraphics[width=0.33\columnwidth]{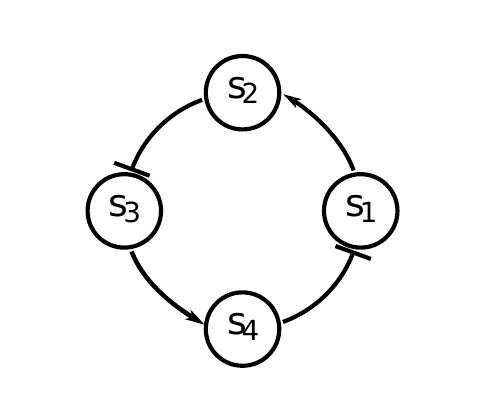} \label{fig:toy-ex}}\qquad
  \subfigure[In this configuration species $S_1$ and $S_3$ are lumped together, while reducing one state]{\includegraphics[width=0.33\columnwidth]{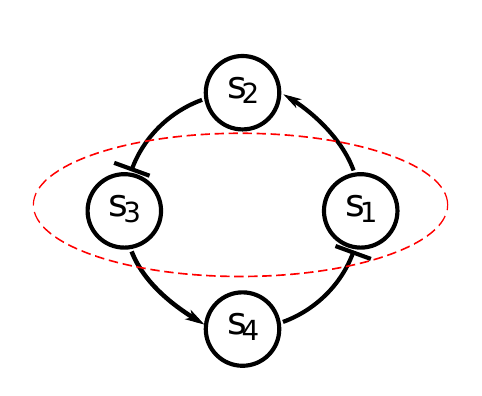} \label{fig:toy-ex-1}}
  \subfigure[In this configuration species $S_1$ and $S_2$ are lumped together, while reducing one state]{\includegraphics[width=0.33\columnwidth]{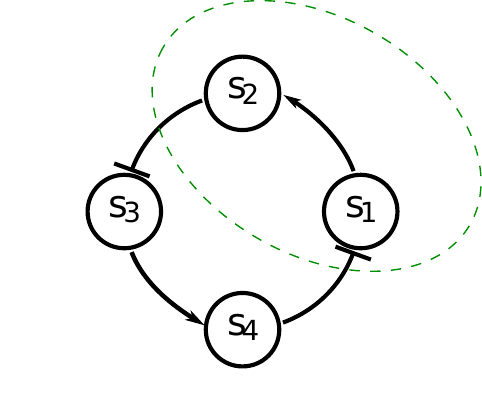} \label{fig:toy-ex-2}}\qquad
  \subfigure[In this configuration the pairs of species $S_1$, $S_2$, and $S_3$, $S_4$ are lumped together, while reducing two states]{\includegraphics[width=0.33\columnwidth]{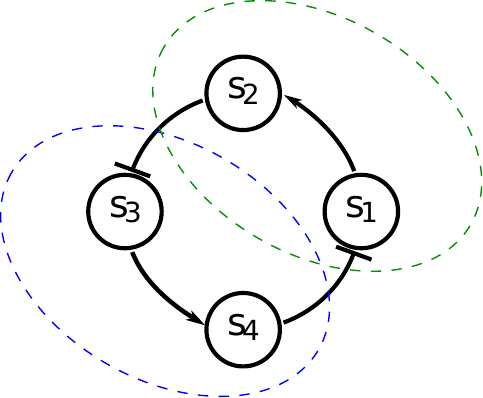} \label{fig:toy-ex-3}}
  \caption{Different configurations for reduction in the toy example.}
\end{figure}
  \begin{align*}
     &\dot m_i = \frac{c_{i 1}}{1 + p_j^2}-c_{i 2} m_i + c_{i 5} u_i\\
     &\dot p_i = c_{i 3} m_i-c_{i 4} p_i
  \end{align*}
where $c_{i1}$ are constants, $m_i$ are mRNA concentrations, $p_i$ are protein concentrations, $u_i$ are exogenous control inputs. $i$ is equal to one or two, $j$ is also equal to one or two, but not equal to $i$. If the state-space is written in the following form $\begin{pmatrix} p_1 & m_1 & p_2 & m_2 \end{pmatrix}$, then this model is monotone with respect to the orthant $\diag(\begin{pmatrix} 1~1~-1~-1\end{pmatrix})\mathbb R^{4}_{\ge 0}$ for all values of parameters. The parameters are chosen as follows:
\[
  c_{1\cdot}  = \begin{pmatrix} 3 & 1 & 1 & 0.2 & 1 \end{pmatrix}\quad c_{2\cdot}  = \begin{pmatrix} 10& 2 & 1 & 0.5 & 1 \end{pmatrix}
\]
This model has two stable steady-states and the state-space is separated into two regions serving as  basins of attraction for these steady-states. We compute the reduced order model using a linearisation around a steady-state $x_{s s} = \begin{pmatrix} 0.14   & 9.8 &   0.03 &   4.9  \end{pmatrix}$, 
and we choose the initial state $x_{0}$ from the basin of attraction of $x_{s s}$:
\[ 
 x_{01}  = \begin{pmatrix} 1  & 10 & 1 & 1 \end{pmatrix}
\]

\begin{table}[t]
  \centering\caption{Reduction of the toy network. The error in the macroscopic concentrations. }
  \label{tab:lin-point}
  \begin{tabular}{cccc}
  Method $\backslash$  Error               & $L_1$ & $L_2$ & $L_{\infty}$\\
\hline
\hline
  QSSA                                      & 67.3 & 11.9 & 3.2\\
  Configuration in Fig.~\ref{fig:toy-ex-1}  & 61.0 & 8.1 & 2.2\\
  Configuration in Fig.~\ref{fig:toy-ex-2}  & 1.9 & 0.59  & 1.1\\
  Configuration in Fig.~\ref{fig:toy-ex-3}  & 13.8 & 2.3  & 0.79\\
\hline
  \end{tabular}
\end{table}

In all the simulations presented in Table~\ref{tab:lin-point}, we set $u=0$, which should give an advantage to the time-scale separation, since in our methods we take into account control signals.
Surprisingly, the difference in the error between QSSA and reduction according to the configuration in \figref{fig:toy-ex-1} is marginal, even though QSSA removes two states and reduction according to the configuration in \figref{fig:toy-ex-1} just one. On the other hand other types of reduction provide much better models if two states (as in the configuration from \figref{fig:toy-ex-3}) or one state (as in the configuration from \figref{fig:toy-ex-2}) are removed. We suppose that the topology of the network has influence on the quality of reduction in this case. The reduction according to the configurations from \figrefs{fig:toy-ex-2}{fig:toy-ex-3} simply removes connections in the network. While the reduction according to the configuration in \figref{fig:toy-ex-1} destroys the topology of the original network.

\subsection{Kinetic Model of Yeast Glycolysis. Non-Monotone Dynamics \label{ex:gly}}
This model was published in~\cite{van2012testing}. It consists of twelve metabolites and four boundary fluxes. In this example, we model the network's response to change of glucose in the system as in~\cite{rao2012model}.  We treat levels of $\mathrm{ATP}$ and glycose $\mathrm{GLCo}$ as control inputs. At time zero we change the levels of $\mathrm{ATP}$ and $\mathrm{GLCo}$ from $3$ to $1.5$ and $0.25$ to $5$ respectively. 

Note that the Jacobian of the dynamics is not Metzler, but there are only five negative off-diagonal elements. Moreover, if we knock out only one one-directional and one bi-directional reaction, then this network will have monotone dynamics with respect to the orthant $\diag(\begin{pmatrix} 1~1~1~1~1~1~1~1~1~1~1~-1\end{pmatrix})\mathbb R^{12}_{\ge 0}$. As was discussed earlier, this phenomenon is not a unique feature of this particular model and it was noticed in \cite{sontag2007monotone}.
Using this intuition, it was not a great surprise that a linearised model around a steady-state would have block-diagonal Gramians with a sparsity pattern according to some state partitioning. However, the existence of diagonal Gramians was a great surprise. This meant that without any reservation we could approximate any group of states, while preserving the other states intact. 

The simulation results are presented in Table~\ref{tab:gly-red} for various reduction configurations. We apply QSSA to metabolite concentrations, while using the lumping method we try to lump those metabolites into one new state, so that the number of reduced states is similar in both cases. First two rows of each subtable in Table~\ref{tab:gly-red} can be compared directly, and it is clear that the proposed reduction method performs better in terms of quality than QSSA. 

The proposed reduction method is also more flexible in terms of reduction choices. In the third row of Subtable~\ref{tab:gly-red}-2, the region \{3PG-PEP\} contains three metabolites; however, we reduced only two states after computing the state-space transformation. In the fourth row, additionally to reducing only one state in region\{3PG-PEP\}, in the region \{GLCi-F6P\} we reduce only one state.  This provides us with the best model among all the reduction attempts. 

\begin{figure}[t]
  \centering
\includegraphics[width=0.66\columnwidth]{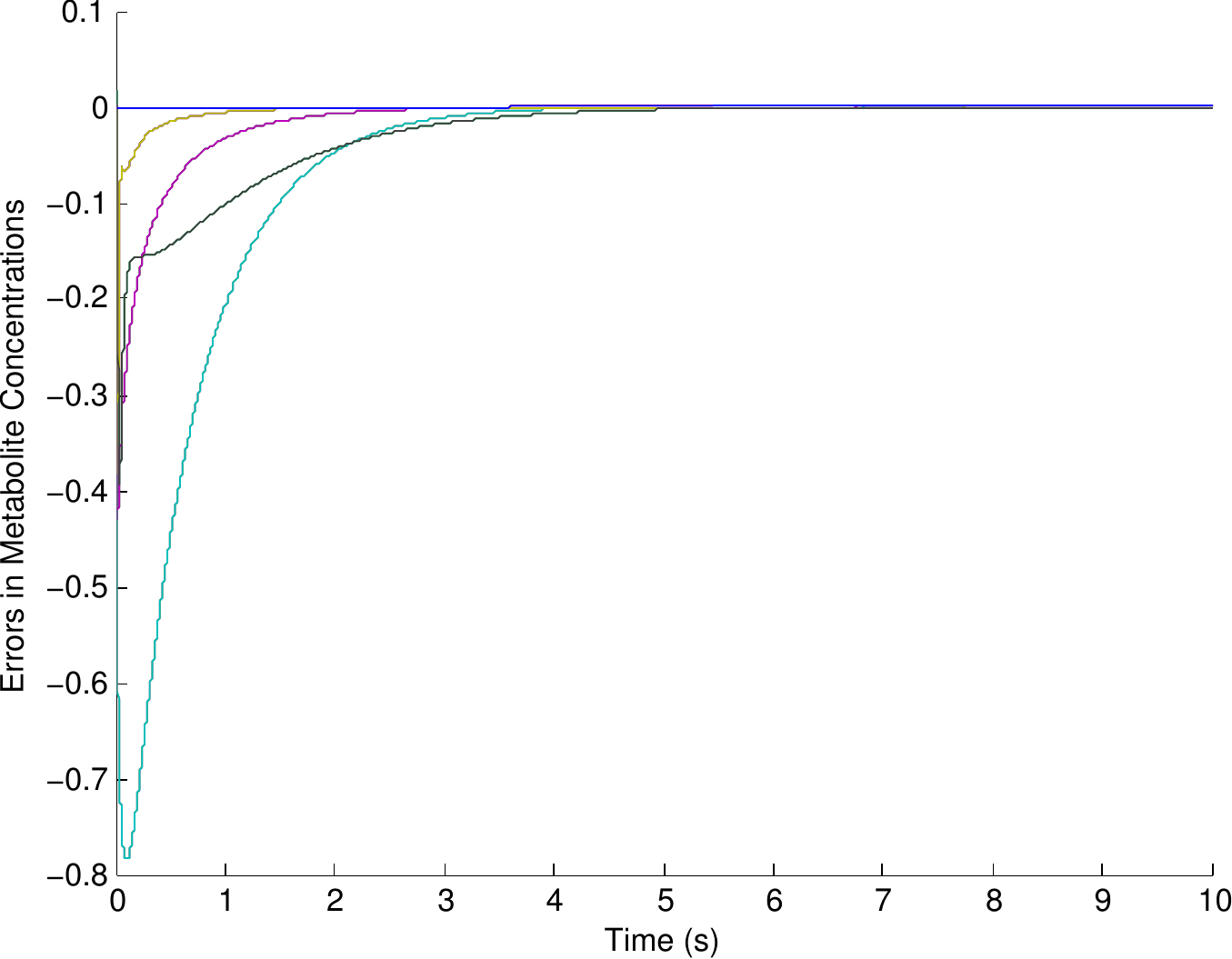}

\includegraphics[width=0.66\columnwidth]{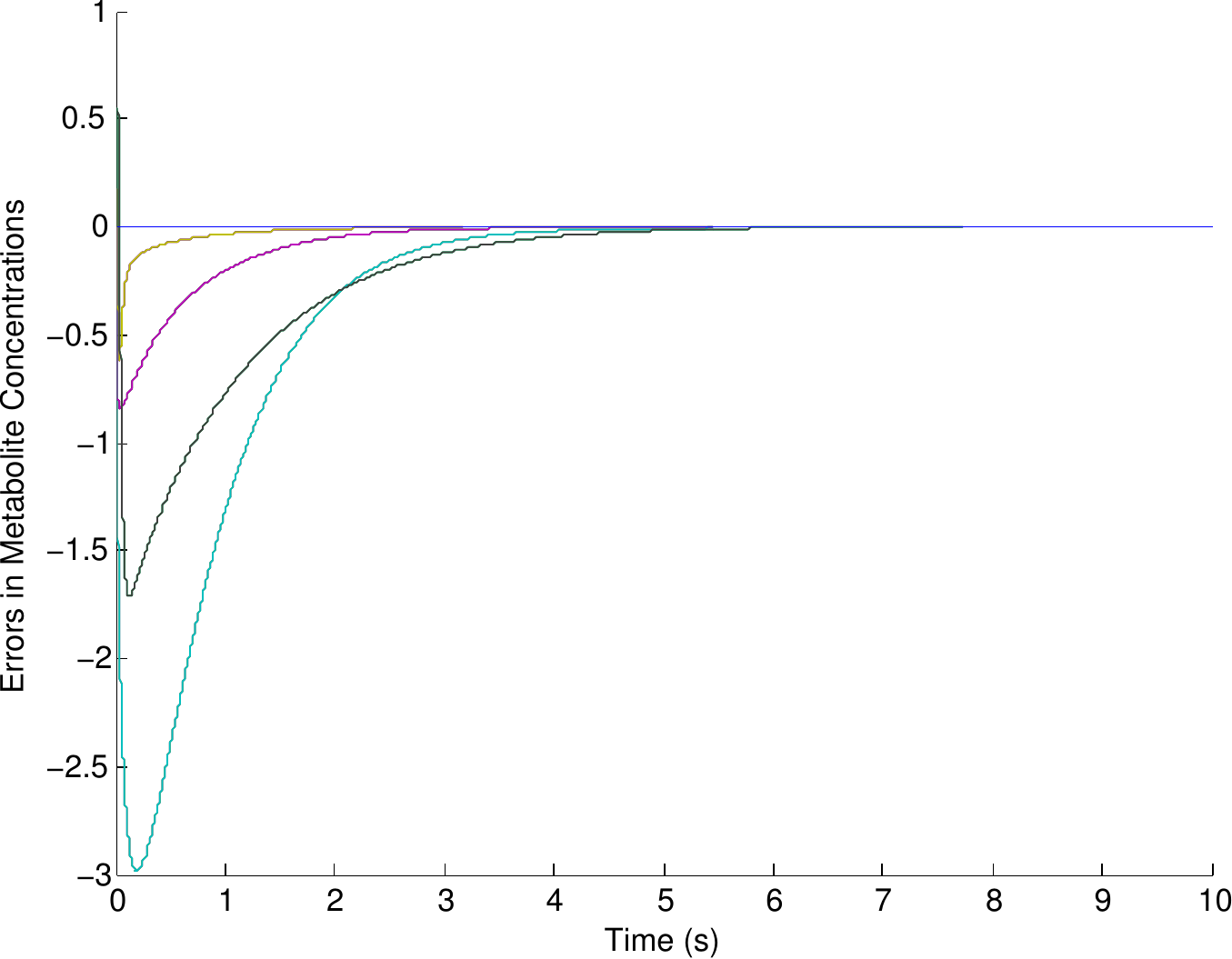}
  \caption{In the upper panel,
the errors between the Method 1 reduced  and the full order models, in the lower panel, the error between the Method 2 reduced and the full order models are depicted. Blue line is the concentration of F16P metabolite, green is TRIO, red - PYR, cyan - AcAld, purple - NADH. 
For the methods 1 and 2 we considered the regions \{GLCi-F6P\}, \{BPG-PEP\} together and reduced two and three states, correspondingly in each region.} 
  \label{fig:gly-red}
\end{figure}

The results of the truncation method (Table~\ref{tab:gly-red}-2) may seem unattractive due to lower approximation quality; however, the difference in terms of qualitative behaviour of the full and the truncated reduced models is not as substantial as the numbers suggest. This is illustrated in \figref{fig:gly-red}. The simulation time of the truncated reduced order model is lower by an order magnitude in comparison with QSSA and the proposed reduction method.

\begin{figure}
  \centering
  \includegraphics[width=0.8\columnwidth]{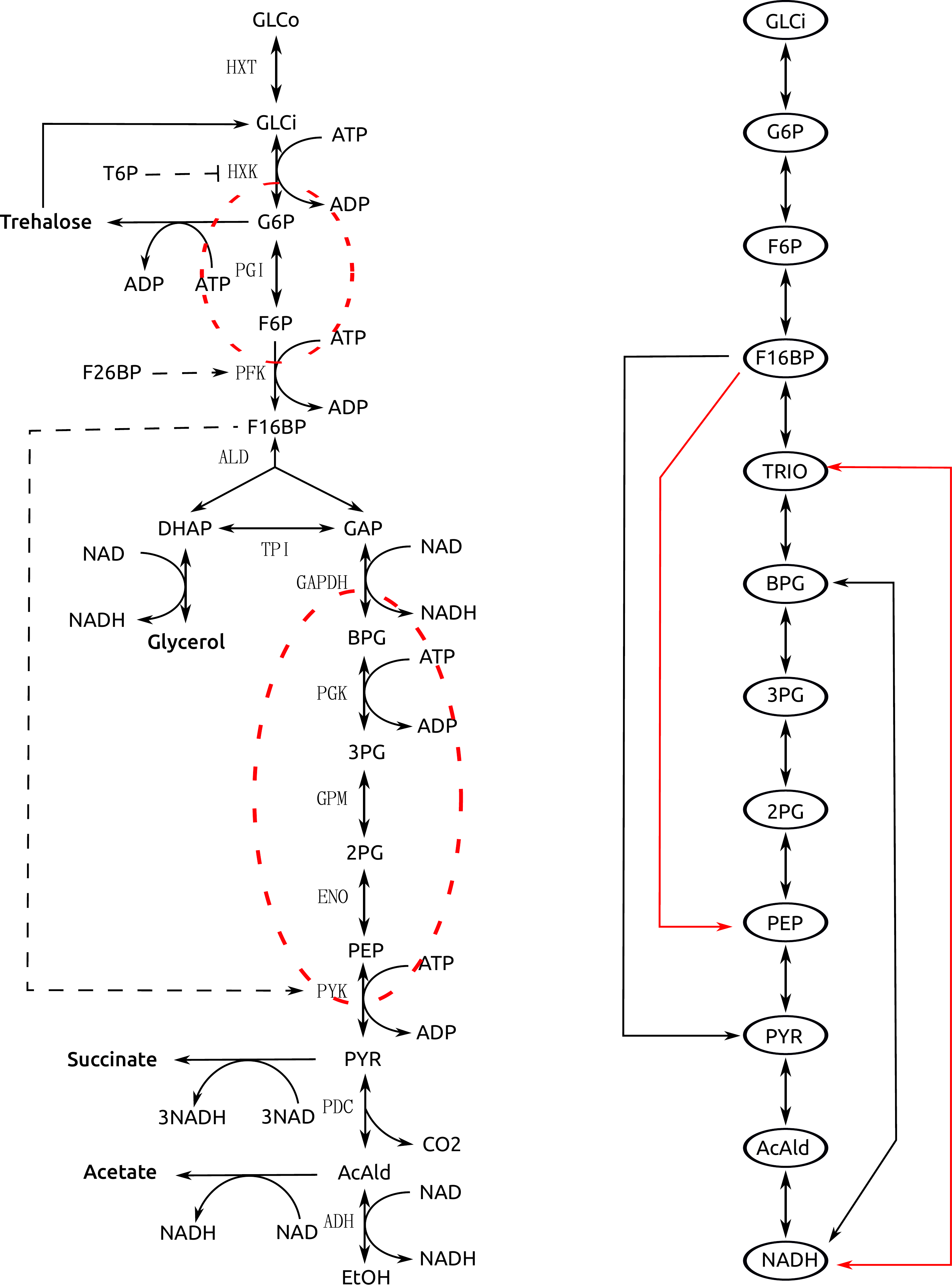} 
  \caption{Depiction of the kinetic model of yeast glycolysis. In the left panel the biochemical graph, and in the right panel a graph of dynamic interactions between metabolites are depicted. If the red connections are removed the dynamics of the network would become monotone.}\label{fig:glycerol}
\end{figure}
 
\begin{table}
\centering\scriptsize
\caption{Deterministic reduction of the glycolysis model. The error of the output is given in different norms. $t$ is simulation time in seconds} \label{tab:gly-red} 
  \begin{tabular}{ccccc}
  \multicolumn{5}{c}{ \sc \ref{tab:gly-red}-1.  QSSA} \\[6pt]
  States $\backslash$ Error & $L_1$ & $L_2$ & $L_{\infty}$& $t$\\
\hline
\hline
  F6P, 2PG, PEP              & $1.21$ & $0.75$ & $0.98$ & $163$\\
  G6P, F6P, 3PG, 2PG, PEP    & $2.05$ & $1.16$ & $1.59$ & $214$\\
\hline
  \end{tabular}
 \vspace{6pt}

  \begin{tabular}{cccccc}
  \multicolumn{6}{c}{ \sc Table~\ref{tab:gly-red}-2.  Reduction by $\{k_1, k_2\}$ states in every region} \\[6pt]
    Lumped Region(s)         & $\{k_1, k_2\}$  & $L_1$  & $L_2$  & $L_{\infty}$ &  $t$\\
  \hline
  \hline
  \{G6P, F6P\}, \{2PG-PEP\} & $\{1, 2\}$ & $1.18$ & $0.79$ & $1.03$ & $161$ \\
  \{GLCi-F6P\}, \{BPG-PEP\} & $\{2, 3\}$ & $1.05$ & $0.57$ & $0.78$ & $260$\\
  \{GLCi-F6P\}, \{3PG-PEP\} & $\{2, 1\}$ & $0.47$ & $0.3$ & $0.4$ & $137$\\
  \{GLCi-F6P\}, \{3PG-PEP\} & $\{1, 1\}$ & $0.14$ & $0.07$ & $0.09$ & $116$\\
\hline\\[-3pt]
  \multicolumn{6}{c}{\sc Table~\ref{tab:gly-red}-3. Truncation by $\{k_1, k_2\}$ states in every region} \\[6pt]
  Lumped Region(s)          & $\{k_1, k_2\}$ & $L_1$ & $L_2$ & $L_{\infty}$&  $t$\\
  \hline 
 \hline
  \{G6P, F6P\}, \{2PG-PEP\}& $\{1, 2\}$ & $15.1$ & $3.2$ & $6.1$ & $14$\\
  \{GLCi-F6P\}, \{BPG-PEP\} & $\{2, 3\}$ & $5.9$  & $2.8$  & $2.9$ & $14$\\
  \{GLCi-F6P\}, \{3PG-PEP\} & $\{2, 1\}$ & $4.1$ & $1.9$ & $1.9$ & $14$\\
  \{GLCi-F6P\}, \{3PG-PEP\} & $\{1, 1\}$ & $4.0$ & $1.8$ & $1.6$ & $15$\\
\hline
  \end{tabular}
\end{table}

\section{Conclusion}
We have presented a method for obtaining structured reduced order models of biochemical reaction networks. The algorithm involves computation of a state-space transformation around a steady-state, followed by a truncation and/or lumping procedure which preserves structure and local monotonicity and stability of the system. The algorithm was illustrated on two numerical examples, one of which was not monotone and compared with a \emph{standard} QSSA based reduction.

\section{Acknowledgment}
The authors would like thank Prof Bayu Jayawardhana and Dr Shodhan Rao for kindly providing the kinetic model of yeast glycolisis. JA acknowledges funding through a junior research fellowship from St. John’s College, Oxford. AS is supported by the EPSRC Science and Innovation Award EP/G036004/1  
\bibliography{Biblio}
\section*{Appendix}
\subsection*{Proof of Lemma \ref{lem:exist}}
It suffices to show that there exist a strictly diagonal $P$ satisfying the controllability Lyapunov inequality as this is a more restrictive case than a block-diagonal and non-negative $P$. Similar arguments hold for a diagonal $Q$ with satisfying the observability Lyapunov inequality. It is known that there exist a diagonal $P$ satisfying the following inequality
\[
    A P + P A^T \le -\delta I
\]
for a positive $\delta$, given an asymptotically stable matrix $A$. Let $X = A P + P A^T$. 
Set $\gamma = {\bar{\sigma}}(B B^T)/\delta $, clearly  $\gamma$ is such that $ \gamma X + B B^T$ is a negative semidefinite matrix. Therefore exist a diagonal $P$ satisfying the Lyapunov inequality, which completes the proof.

\subsection*{Proof of Lemma \ref{lem:main}}
\begin{enumerate}
\item The existence of a balancing transformation is an established result (cf. \cite{RCBZhou}). $P_{22} Q_{22}$ is an irreducible matrix with nonnegative entries, therefore by Perron-Frobenius theorem there exist a positive eigenvector $w$ such that 
\[
P_{22} Q_{22} w = \Sigma_{11}^2 w
\]
where $\Sigma_{11}^2$ is the entry $(1,1)$ of the matrix $\Sigma^2$ and the largest eigenvalue of $P_{22} Q_{22}$. From the existence of $T$, it follows that $T \Sigma^2 T^{-1} = P_{22} Q_{22}$, where $\Sigma^2$ is a diagonal matrix. Hence $T$ contains right eigenvectors to a matrix $P_{22} Q_{22}$ and without loss of generality $w$ is the first column of $T$.  Similarly it can be shown that $v$ the first column of $T^{-T}$ is nonnegative. 
\item Stability of the matrix $A_t$ is a collection of known results, but it is presented for completeness. Let $T$ be $\blkdiag(I_{n-k},T_{22})$. Introduce the following partitioning of these matrices:
\end{enumerate}
\[
\begin{gathered}
 T^{-1} A T =  \begin{pmatrix}  A_t & A_{tr} \\ A_{rt} & A_{rr} \end{pmatrix} \quad 
T^{-1} P T^{-1} = \begin{pmatrix} \tilde P & 0 \\ 0 & \tilde \Sigma \end{pmatrix}\\
\tilde P  = \begin{pmatrix} P_{11} & 0 \\ 0 & \Sigma_{1 1} \end{pmatrix}\quad
\tilde \Sigma = \diag\left(\Sigma_{2 2}, \Sigma_{n n}\right)  
\end{gathered}
\]
Now stability of $A_t$ can be established by simply writing the Lyapunov inequalities in the new variables.
\begin{equation*}
  \begin{aligned}
  T^{-1} A T T^{-1} P T^{-T} + T^{-1} P  T^{-T} T'  A' T' &\le 0 \\ 
    \begin{pmatrix} A_{t} & A_{tr} \\ A_{rt} & A_{rr} \end{pmatrix} \begin{pmatrix} \tilde P & 0 \\ 0 & \tilde \Sigma \end{pmatrix}   + \begin{pmatrix} \tilde P & 0 \\ 0 & \tilde \Sigma \end{pmatrix} \begin{pmatrix} A_{t}' & A_{rt}'  \\ A_{tr}'  & A_{rr}' \end{pmatrix}  &\le 0 \\
  \begin{pmatrix} A_{t} \tilde P + \tilde P A_{t}' & A_{tr} \tilde \Sigma + \tilde P A_{rt}'\\ \ast &
 A_{rr} \Sigma + \Sigma A_{rr}' \end{pmatrix} &\le 0
  \end{aligned}
\end{equation*}
Proving that $A_t$ is Metzler is also straightforward. $A_{12} w$ $v' A_{21}$ are nonnegative since $w$, $v$, $A_{12}$, $A_{21}$ are individually nonnegative \cite{grussler2012symmetry}.  All is left to show that $v' A_{22} w$ is a negative scalar. Since 
$A_{t} \tilde P + \tilde P A_{t}' \le 0$ then
\[
\begin{pmatrix} A_{11} P_{11} + P_{11} A_{11}' & A_{12} w \Sigma_{11} + P_{11} A_{21}' v\\ \ast &
 v' A_{22} w \Sigma_{11} + \Sigma_{11} w' A_{22}' v \end{pmatrix} \le 0
\]
and hence $v' A_{22} w \Sigma_{11} + \Sigma_{11} w' A_{22}' v$ is negative, which implies that $v' A_{22} w$ is negative since $\Sigma_{11}$ is a positive number.
\begin{enumerate}  \setcounter{enumi}{2}
\item This result is shown in~\cite{Sandberg09}.
\end{enumerate}
  \end{document}